\documentclass[11pt]{article}
\usepackage{amsbsy,amsmath,amsthm,amsfonts,amssymb,amstext,amscd,longtable,flafter,float}

\newtheorem{theo}{\bf Theorem}[section]

\newtheorem{defi}{\bf Definition}[section]
\newtheorem{nota}{\bf Notation}[section]
\newtheorem{lem}{\bf Lemma}[section]

\begin{document}

\begin{center}
{\bf {\Large  Optimality of multi-way designs}}

\vskip5pt
{\bf email of corresponding author : bhaskarbagchi53@gmail.com}

\vskip5pt

{\bf {\large  Sunanda Bagchi and  Bhaskar Bagchi*\\[0pt]
1363, 10th cross, Kengeri Satellite Town, \\
Bangalore 560060 \\
India\\}}
\end{center}

\vskip10pt

 {\bf {\large Abstract  }}

\vskip5pt

In this paper we study optimality aspects of  a certain type of designs in a multi-way  heterogeneity setting.  These are ``duals" of
plans orthogonal through the block factor (POTB). Here by the dual of a main effect plan (say $\rho$) we mean a design in a  multi-way
heterogeneity setting obtained from $\rho$ by interchanging the roles of the block factors and the
treatment factors. Specifically, we take up two series of universally optimal POTBs  for symmetrical experiments
constructed in Morgan and Uddin (1996). We show that the duals of these plans, as multi-way designs, satisfy M-optimality.

Next, we construct another series of multiway designs and proved their M-optimality, thereby generalising the result of Bagchi and Shah
(1989). It may be noted that M-optimality includes all commonly used optimality criteria like A-, D- and E-optimality.

\vskip5pt


  %Submission information :    bhaskarbagchi53@gmail.com
% JSTP : username : BhaskarBagchi, Pass word : Bachchu1953

\section{Introduction}
An experimental unit  subjected to more than one heterogeneity directions occurs in many situations. The search for optimal designs in a
setting with two crossed nuisance (or blocking)  factors began in Kiefer (1975). Later others joined the search and obtained optimal
row-column designs [see Shah and Sinha (1989) for more details].

Optimality study in a general $m$-way heterogeneity setting was initiated by Cheng (1978). He assumed a model with no interaction
    among the  blocking factors and  considered a set up with constant number of level combinations of the  blocking factors. Mukhopadhyay
     and  Mukhopadhyay (1984) assumed  the same model  as Cheng (1978), but relaxed  the complete crossing requirement on the block factors.
    Assuming that the level combinations of the block factors form an  orthogonal array with variable number of symbols, they obtained
    optimality results very similar to those of Cheng (1978). Bagchi and  Mukhopadhyay (1989) considered the situations where two factor
    interactions are present among the block factors and obtained optimality results. Morgan (1997) considered an m-way setting with
    $t$-factor interactions and generalised all the optimality results obtained so far.

     The dual of an optimal  block design often satisfies optimality property. In this paper we try to see whether
    that happens for a main effect plan (MEP) also. In this context we interpret duality as an interchange between the roles of the block factors and the treatment factors. Thus, dual of a blocked MEP is a design in a  multi-way setting  [see Definition \ref {dual}]. We have taken up two series of universally optimal blocked MEPs
    constructed in Morgan and Uddin (1996) and proved that the duals of these plans are also optimal as multi-way  designs, although the
    type of optimality is different and the competing class of designs is smaller.
          We have also constructed a new series of  multi-way  designs and proved its optimality.

 In section 2 we describe the set up we consider. In this section we also discuss the concept of adjusted orthogonality and its consequence [see Definition \ref {adjustedOrth} and Lemma \ref {adjOrth}],  which we used later. Discussion on optimality criteria is placed in Section 3.
  In Section 4 we describe three multiway settings and derive the information matrix (or C-matrix) of the treatment effects for a design in each of these settings [see Lemma \ref {C-matTreatment}]. We take up the designs of our interest in Section 5. We describe the duals of
  two series of plans of Morgan and Uddin (1996), study the properties of these multiway designs and then prove their optimality [see Theorems
\ref {opt d1*} and \ref {d2*}]. In Section 5.3 we construct a series of multiway designs satisfying adjusted orthogonality [see Theorem \ref
{propOfd^*3}] and proved their optimality property [see Theorem \ref {optd3*}].

\section{Preliminaries}
%\subsection{The set up of a multifactor experiment}

We shall consider designs in a  multi-way heterogeneity setting, in which an experimental unit is subjected to more than one heterogeneity
directions or blocking factors.
 An experimental design in a  multi-way heterogeneity setting is completely specified by a quadruple
$(U, {\cal F}, S, \eta)$, where $U$ is the set of all units, $ {\cal F} = {\cal B} \cup \{V,G\}$ is the set of  factors. The distinguished
factor $V$ is called the treatment factor, the factors in ${\cal B}$ are called the block factors and $G$ is the general effect. $S$
represents the set of levels : $S_A$ being  the set of levels of the factor $A, \; A \in {\cal F}$. $s_A : = |S_A|$. In particular, $S_V$ is
the set of treatments. We shall write $v$ for $s_V$. We have $s_G = 1$.

The allotment function $\eta : U \times  {\cal F} \rightarrow S$ specifies, for any unit $u$ and factor $A$, the level $\eta(u,A)$ of  $A$ to be applied to the unit $u$.
Thus, $\{\eta(u,A) : u \in U\} = S_A, \; A \in {\cal F}$.

 $m$ will denote  $ |{\cal B}|$ and $n$  will denote $|U|$. We assume that $m \geq 3$.
 %[A set up with $m = 2$ is commonly referred to as a row-column set up.]

  \begin{nota} \label{allotment}  $\eta_1$ (respectively  $\eta_2$) will denote the restriction of the allotment function $\eta$ to $U
  \times \{V\}$ (respectively to $U \times  {\cal B}$).  $\eta_2$ is called the experimental set up (and is viewed as given by nature) and
  $\eta_1$ is called the design (and is viewed as designed by the experimenter). Note that  $\eta_1$ may be thought of as the function $x
  \mapsto \eta_1(x,V)$ from $U$ to $S_V$. For ease of description, we shall describe a design by the allotment function $\eta$, which
  includes both the experimental set up and the allocation of treatments in a combined form.
  \end{nota}

 \begin{nota} \label{model}
 (a)   $ 1_n$ will denote the $n \times 1$ vector of all-ones, while $
J_{m \times n }$ will denote the $m \times n$ matrix of all-ones. We shall write $J_m$ for $J_{m \times m}$.
 The $n \times 1 $ vector of responses will be denoted by $Y$.
%[That is  the number of blocking factors is denoted by $m$ and the number of experimental units by $n$]

(b)  Fix $A \in  {\cal F}$. The design matrix $X_A$ is the $n  \times s_A$   matrix having the $(u,p)$th entry  $1$ if
 $\eta(u,A) = p$,  and $0$ otherwise ($ u \in U, p \in S_A$).  Clearly, $X_G = 1_n$.
We shall use the following notations for the incidence matrices.  For $A \neq B$ in   ${\cal B}$, the $A$ versus $B$ incidence matrix is
the matrix  $M_{AB} = X'_A X_B$. Thus the $(p,q)$th entry of $M_{AB}$ is  $M_{AB} (p, q)$, which is the number of units $u$ in which $A$ is
at level $p$ and $B$ is at level $q, \;p \in S_A, \; q \in S_{B}$.  $A$ is said to be equireplicate if the levels of $A$ appear with equal
frequency in the setting.

  For $A \in {\cal B}$, the $V$ versus $A$ incidence matrix is  the matrix  $N_{A} = X'_V X_A$.

  (c) The  $ s_A \times 1$ vector $\alpha^A$ will denote the vector of unknown effects of $A, \;A \in {\cal F}$. \end{nota}
%The scalar $\alpha^0$ will denote general effect.

 The model is expressed in matrix form as
\begin{equation}\label{modelEq} {\mathbf Y} =  {\mathbf X} \alpha +\epsilon :
%\mbox{ where }
{\mathbf X} = \left[ X_A : \; A \in {\cal F}\right] \;\: \alpha =
\left[ \alpha^A : \; A \in {\cal F } \right ]^\prime, \;\;  \epsilon \sim N_n(0,\sigma^2 I_n).
\end{equation}
Here $Z \sim N_p(\mu,\Sigma)$ means that $Z$ is a random variable
following $p$-variate normal distribution with mean $\mu$ and
covariance matrix $\Sigma$.

  The normal equation for the least square estimates of the vector $\alpha$ of all effects is
 \begin{equation}  \label{NEbeta}
  {\mathbf X}'  {\mathbf X} \; \widehat{\alpha} = {\mathbf X}' {\mathbf Y}.
  \end{equation}

We shall use the following notations for the sake of compactness.
\begin{nota} \label{C-matandQ}
(a) For any  matrix $M$,  $P_M$ will denote the projection operator onto the column space of $M$. Thus, $P_M =M (M'M)^-M'$, where $H^-$
denotes a g-inverse of the matrix $H$.

(b)  Let   $T = \{A,B, \cdots\}$ be a  subset of ${\cal F}$.

 (i) $ X_T$ will denote $ \left[ \begin{array}{ccc}{\mathbf X}_A &
 {\mathbf X}_B & \cdots \end{array} \right ]$.

 %Moreover, $X'_A X_T$ will be denoted by $M_{AT}$.

  (ii) $\alpha^T$ will denote  $ \left[ \alpha^A)' : A \in T \right ]^\prime$. $\widehat{\alpha^T}$ will denote the least square
estimate of $\alpha^T$.

(iii)  $ P_A$ will denote the projection operator onto the column
space of $ X_A, A \in {\cal F}$. Further, $ P_T$ will denote the projection
operator onto the column space of $ X_T$.
\end{nota}

Consider  $ T \subset {\cal F}$  and let $\bar{T} = {\cal F} \setminus T$.
The reduced normal equation for the vector $\alpha^T$ of  all effects of all members in $T$ after eliminating the effects of $\bar{T}$
from (\ref {NEbeta}) is given by
\begin{equation}  \label{ReducedNEalphaT0}
 {\mathbf C}_{T;\bar{T}} \widehat{\alpha^T} = {\mathbf Q}_{T;\bar{T}}, \mbox{ where }
 \end{equation}
\begin{eqnarray} \label{CQSS}
 {\mathbf C}_{T;\bar{T}} & = & (( C_{AB;\bar{T}} ))_{A,B \in T},\;  C_{AB;\bar{T}}  =  {\mathbf X}'_A (I - P_{\bar{T}}) {\mathbf X}_B, \\
{\mathbf Q}_{T;\bar{T}} & = & (( Q_{A;\bar{T}} ))_{A \in T}, \;  Q_{A;\bar{T}} = {\mathbf X}'_A (I - P_{\bar{T}}) {\mathbf
Y}.\end{eqnarray}

\begin{defi}\label{CdAdjTotal}
(a) For factors $A \neq B$,  $Q_{A;B}$ will denote the vector $ X_A^\prime (I-P_B)Y$. It is commonly known as ``the total for $A$, adjusted for $B$".

  (b) By the C-matrix of a multi-way design $d$ we mean the $v \times v$ matrix ${\mathbf C}_{V;\bar{V}} = {\mathbf X}'_V (I - P_{\bar{V}})
  {\mathbf X}_V$, which we shall refer to as $C_d$. In order that every treatment contrast  is estimable, rank of $ C_d$ must be $v-1$. We,
  therefore, consider only the designs satisfying that condition. We shall refer to such designs as `connected'.
\end{defi}

%For ready reference we expand the expression for $C_{AB;G}$ below.
%\begin{equation}  \label{EliminatedG} C_{AB;G} = M_{AB} - (1/n) r_A' r_B. \end{equation}

\begin{defi}\label{equireplicate} A  multi-way design $d$ is said to be equireplicate if every treatment (i.e. every member of $S_V$)
occurs with equal frequency.
\end{defi}

{\bf Adjusted orthogonality :} Following Eccleston and Russel (1977) we define the following.
\begin{defi}\label{adjustedOrth} For three distinct factors $A,B,T$ we say that $A$ and $B$ are  adjusted orthogonal with respect to $T$ if
$Cov(Q_{A;T}, Q_{B;T}) =0$.

If  $T$ is  equireplicate, then $Cov(Q_{A;T}, Q_{B;T}) = M_{AB} - r^{-1} M_{AT} M_{TB} $, where $r$ is the replication number of $T$.
Therefore,in this case, $A$ and $B$ are  adjusted orthogonal with respect to $T$, if and only if $M_{AT}M_{TB} = r M_{AB}$.
\end{defi}

 Shah and Eccleston (1986) proved a few interesting properties of a row-column design, in which the factors row and column are
 adjusted orthogonal with respect to the treatment factor. Those results can be easily generalized to a multi-way  heterogeneity setting.
 Among those, we present one result, restricting to the designs in the equireplicate class.

 \begin{lem}\label{adjOrth} Consider an equireplicate multi-way design $d$. Suppose $ A,B \in {\cal B}$ are such that $M_{AB} =
 J_{s_A\times s_B}$. Consider the $v \times v$ matrices $T_A = N_{A}N^\prime_{A}$ and $T_B = N_{B}N^\prime_{B}$. If $A$ and $B$ are
 adjusted orthogonal with respect to $V$, then the following holds.

Suppose $x'1_v =0$ and $x$ is an eigenvector of $T_A$ with non-zero eigenvalue. Then $T_B x = 0$.
 \end{lem}

 {\bf Proof :}   Since  $d$ is equireplicate, $1_v$ is an eigenvector of $T_A$ as well as of $T_B$  corresponding to their maximum eigenvalue.  Again, by the hypothesis,  $T_A T_B = r^3 J_v = T_B T_A$, where $r$ is the replication number (for the treatments). Thus, $T_A$ and $T_B$ are commuting matrices and hence there is an orthonormal basis (including $1_v$) consisting of common eigenvectors of these two matrices.  Therefore,
       ${\cal C} (T_A) \cap {\cal C} (T_B) = {\cal C} (T_A T_B) = {\cal C} (J_v) = \langle 1_v \rangle$. Hence the result. $\Box$

\section{Optimality Tools}
\setcounter{equation}{0}
We define a class of  optimality criteria which are functions of the eigenvalues of the C-matrix of a design $d$ [see Definition \ref
{CdAdjTotal}].

\begin{nota}
\label{Eigenvalue} For a real symmetric $n \times n$ matrix $A$,

(a) $ \mu_0 (A) \leq \cdots \leq \mu_{n-1} (A)$ will denote the
eigenvalues of $A$.

(b) $\mu(A) = (\mu_0 (A), \cdots \mu_{n-1} (A))^{\prime }$ will
denote the vector of eigenvalues of $A$.

(c) Consider a connected design $d$. $C_d$ will denote its C-matrix. By the vector of eigenvalues of $d$ we shall mean the vector of
positive eigenvalues of $C_d$ in the nondecreasing order and it will be denoted by $\mu(d)$.
\end{nota}

\begin{nota}\label{Psif}
  $\Phi$ will denote the class of all  non-increasing  convex real valued functions on $I\!\! R^+ = (0, \infty)$.

For $f \in \Phi, \: \Psi_f : (I\!\! R^+)^n\rightarrow
I\!\!R$ is defined by $\Psi_f (x) = \sum_{i=1}^n f( x_i)$.
\end{nota}

\begin{defi} \label{PsiOpt} (a) For   $x,y \in (I\!\! R)^n$,  $x$ is said to be $\Psi_f$-better than $y$ if
\begin{equation}  \label{PsiOptEq} \Psi_f (x) \leq \Psi_f (y).
\end{equation}

(b) Consider a class $\mathcal{D}$ of connected designs with a common set up and a member $d^*$ of $\mathcal{D}$.  If $\mu(d^*)$ is
$\Psi_f$-better than $\mu(d)$ for every $d \in \mathcal{D}$, then $d^*$ is said to be $\Psi_f$-optimal in $\mathcal{D}$.
\end{defi}

The following members of $\Phi$ are  of special interest as the
corresponding $\Psi_f$ criteria have important statistical interpretation [see Shah
and Sinha (1989) for more details]. These are
the functions $f(u) \equiv u^{-1}$ corresponding to A-optimality and
$f(u) \equiv -log (u)$  corresponding to D-optimality. Another popular
optimality criterion is E-optimality, which may be obtained as the limit of
a class of $\Psi_f$ criteria.

A powerful approach to design optimality problems
% when it applies, encompasses all members of ${\cal F}$
is through the concept of majorization.

%$(I\!\! R^n)^{\uparrow}$

\begin{nota}\label{xUpwards} $(I\!\! R^n)^{\uparrow} $ will denote the set of vectors  $x = (x_1, \cdots x_n) \in I\!\! R^n$,
such that $x_1 \leq \cdots \leq x_n$.
 \end{nota}

\begin{defi} (Marshall, Olkin and Arnold(2011))
\label{weakMajor} For $x,y \in (I\!\! R^{n})^\uparrow $, $x$ is said to be
weakly majorized from above by $y $ (in symbols, \textbf{$x \prec^w
y$}) if
\begin{equation}  \label{major}
\sum_{i=1}^{k} x_{i} \geq \sum_{i=1}^{k} y_{i}, ~~ k =
1,2, \cdots, n, %\; \mbox{ with `$=$' for } k = n.
\end{equation}
\end{defi}

\noindent
See Marshall, Olkin and Arnold (2011) for a comprehensive treatment of majorization concepts and results.

Following Bagchi and Bagchi (2001) we define the following.
\begin{defi} \label{M-opt}
 For $d, d^{\prime } \in \mathcal{D}$, $d$ is said to be better than $d^{\prime }$ in the sense of majorization (in
short M-better),  if $\mu(C_d)  \prec^w \mu(C_{d^{\prime }})$, equivalently if $\mu(d) \prec^w \mu(d^{\prime })$.  A design $d^*$ is  said
to be optimal in $\mathcal{D}$ in the sense of majorization
%in a subclass $\mathcal{D}$ of $\mathcal{D}_{b,k,v}$
(or, in short, $d^*$ is M-optimal in $\mathcal{D}$) if it is M-better
than every member of $\mathcal{D}$.
\end{defi}

In view of a theorem in Tomic (1949), we have the following useful result. [See Proposition B.2 of chapter 4 of
Marshall, Olkin and Arnold (2011) for the theorem of Tomic ].

\begin{theo} \label{M-optCond}  A design $d^* \in \mathcal{D}$ is M-optimal in $\mathcal{D}$ if and only if $d^*$ is
 $\Psi_f$-optimal in $\mathcal{D}$ for every $f \in \Phi$.
  \end{theo}

By the discussion after Definition \ref  {PsiOpt}, M-optimality implies A-, D- and E-optimality, among many other optimality criteria.

 {\bf Sufficient conditions for  majorization :}. It is easy to verify the following.
\begin{lem}\label{suffMopt0} Suppose $A$ and $B$ are real symmetric matrices. If $ A \geq B$,
 then  $\mu (A) \prec^w \mu(B)$. \end{lem}

The following is another useful result.
\begin{lem}\label{suffMopt} Let $m < n$. Consider  $x \in (I\!\! R^n)^{\uparrow} $, such that $m$ of the  $x_i$'s
are $= a$ and the remaining $n - m$ of them are $=b$, where $a < b$. Suppose $y \in (I\!\! R^n)^{\uparrow} $ satisfies $(1) :
\sum_{i=1}^{n} y_i  \leq  \sum_{i=1}^{n} x_i$ and $(2) :  y_{m+1} \geq b$.  Then,  $x \prec^w y$.
%each of the following conditions is sufficient for \begin{equation}\label{yProp}    \mbox{ or }  y_{m}  \leq  a.\end{equation}
  \end{lem}

  {\bf Proof :} We shall use the fact that for an $ n\times 1$ vector $y \in (I\!\! R^{n})^\uparrow,\; \; (3) : (1/i) \sum_{j=1}^{i} y_{j}
  \mbox{ is increasing in } i. $

  (1) and (2)  together  implies that  $\sum_{j=1}^{m} y_j \leq ma$. This, in view of (3) implies that $\sum_{j=1}^{i} y_j \leq ia, \mbox{
  for } 1 \leq i \leq m.$

  Next, let $m + 1 \leq i \leq n$. Since $y \in (I\!\! R^n)^{\uparrow}$, (2) implies $\sum_{j=i+1}^{n} y_j \geq (n-i)b$. Subtracting this
  inequality from (1), we get $\sum_{j=1}^{i} y_j \leq ma + (i - m)b$. $\Box$

\section{C-matrices of  multi-way designs in a nearly orthogonal set up}
\setcounter{equation}{0}

Our aim is to study the performances of a few multi-way designs. For that we need the  C-matrix of the reduced normal equations for the
BLUEs of the treatment effects for designs in the appropriate  settings, (which is ${\mathbf C}_{V;\bar{V}}$ in Definition \ref {CdAdjTotal} (b)).
 The C-matrix of a design in a multi-way set up has been derived in Mukhopadhyay and Mukhopadhyay (1981), under the assumption that  blocking
 factors are orthogonal to each other. This  condition is not
satisfied in any of the experimental set ups  we consider here and so those results are not applicable.  However, they satisfy certain
other conditions, which may be described as `near orthogonality', as the incidence matrix of every pair of blocking factors is of the form
$aI + bJ$. Here $a$ and $b$ are
integers varying from one type of setting to another.

We now describe  three different types of settings. In this description, the number of blocking factors is assumed to be an arbitrary
integer $m \geq 3$.
%In each of these settings, each level of each blocking factor appears equally often.
\begin{nota}\label{3setUp} (a) By setting of type 1 we mean a  setting with  $n = s(s-1)$ units, in which each blocking factor has s
levels. Moreover, $M_{BB'} = J_s-I_s$ for any two distinct block factors $B, B'$.

(b) A  setting of type 2 has $n = s(p+s)$ for an integer $p$, $s$ levels for each  blocking factor  and $M_{BB'} = p I_s +  J_s, \;  B'
\neq B, \in {\cal B}$.

(c) A  setting of type 3 has $n = s(s+1)$, one (blocking) factor (say $E$) having $s+1$ levels, while the other ones have $s$ levels each.
Let $\tilde{{\cal B}} =
{\cal B} \setminus \{E\}$. The incidence matrices are as follows.
 \begin{eqnarray} \label{incMatType3} M_{BB'} = I_s + J_s, & B' \neq B \in \tilde{{\cal B}},\\
                                    M_{BE} = J_{s \times s+1}, & B \in \tilde{{\cal B}},.\end{eqnarray}
\end{nota}

We proceed to obtain the reduced normal equation for the treatment effects. For the sake of compactness, we shall use the following
notations in the statement of Lemma \ref {C-matTreatment} (the original notations of Section 2 will be used whenever necessary).
\begin{nota}\label{setUpProp} (a) $C_0 = C_{VV;G}$ [see (\ref {CQSS})]. Thus, $C_0$ is the $v \times v$ matrix
$R - (1/n) {\mathbf r} {\mathbf r}'$. Here ${\mathbf r} = (r_1, \cdots r_v)'$ is the vector of replication numbers of the treatments and $R$ is the diagonal  matrix with  diagonal entries same as the entries of ${\mathbf r}$ in the same order.

 (b) $C_B = C_{BV;G}, \; B \in {\cal B}$.
 % [see the equation next to (\ref {ReducedNEalphaT})].

% (c) Finally, $S_{{\cal B}}$ will denote $ \sum \limits_{B \in {\cal B}} C_B$.
 % $C_B$ as in (b) above.
 \end{nota}
 We note that in each of the above settings, ${\mathbf X}'_B 1_n = (n/s_B) 1_{s_B}, \forall B \in  {\cal B}$. Hence
 \begin{equation}\label{CBcolsum0} C_B = N_B - (1/s_B){\mathbf r} 1'_{s_B} \mbox{ and so  } 1'_v C_B = 0.\end{equation}

\begin{lem}\label{C-matTreatment} In the multi-way settings of the three types described in Notation \ref {3setUp} the  C-matrix
 of a design $d$  is as given below.

(a) If  the setting is of type 1 or type 2, then the  C-matrix  is
%$$C_d = C_0 - (1/s) \sum \limits_{B \in {\cal B}} C_B C'_B- (1/(s(s-m)) S_{{\cal B}} S'_{{\cal B}}. $$

%(\sum_{i=1}^{m} N_i)(\sum_{i=1}^{m} N'_i) \vskip5pt

%(b) If the setting is of type 2, then, $C_d$ is as follows.

$$C_d = C_0 - \frac{1}{s} \sum \limits_{B \in {\cal B}} C_B C'_B + \frac{p}{su} S_{{\cal B}} S'_{{\cal B}}, \mbox{ where } p = -1 \mbox{ for Type
1 }, u = s + mp, \mbox{ and } S_{{\cal B}} = \sum \limits_{B \in {\cal B}} C_B. $$

%Check  the expression.

(b) The C-matrix for a design in a setting is of type 3 is as given below.
$$C_d = C_0 - \frac{1}{s} \sum \limits_{B \in {\cal B}} C_B C'_B  + \frac{1}{su} \tilde{ S_{{\cal B}}} \tilde{ S_{{\cal B}}}' ,$$
 where $u = s + (m - 1) $ and $ \tilde{ S_{{\cal B}}} =  \sum \limits_{B \in \tilde{{\cal B}}} C_B, \; \tilde{{\cal B}}$ as in Notation
 \ref
 {3setUp} (c).

 %Check  the expression.
\end{lem}

{\bf Proof :} Taking $T = {\cal F} \setminus G$ in (\ref {ReducedNEalphaT0})  we get the following system of reduced normal equations
(written in compact form).
\begin{eqnarray} \label{NEeliminate0} C_{{\cal B};G} \widehat{\alpha^{{\cal B}}} + C_{{\cal B}V;G} \widehat{\alpha^V}  &=& Q_{{\cal B};G}\\
\mbox{ and } C_{V{\cal B};G}  \widehat{\alpha^{{\cal B}}} +  C_0 \widehat{\alpha^V} &=& Q_{V;G}. \end{eqnarray}
Here $ C_{{\cal B};G} = (( C_{AB;G}))_{A,B \in {\cal B}}$, $C_{{\cal B}V;G}$ is the $ms \times v$ matrix $= [ C_{BV;G} :B \in {\cal
B}]'$ and $C_{V{\cal B};G}$ is the transpose of $C_{{\cal B}V;G}$.

Note that (\ref {NEeliminate0}) consists of $m$ equations, a typical  one of which is
\begin{equation}\label{ith} C_{BB;G} \widehat{\alpha^B} + \sum \limits_{B' \neq B} C_{BB';G} \widehat{\alpha^{B'}}
+ C_{BV;G} \widehat{\alpha^V} = Q_{B;G}. \end{equation}

\begin{equation} \label{Kmatrix}\mbox{ For an integer }t \mbox{ let } K_t \mbox{ denote the symmetric idempotent matrix } I_t - (1/t) J_t .
\end{equation}
{\bf Proof of (a):} By the hypothesis,  we find that the individual C-matrices involving only the blocking factors are as follows.
\begin{equation}\label{CmatricesIn(b)} C_{BB';G}  =  \left\{ \begin{array}{ll}
(p+s)K_s & \mbox{ if }  B' = B, \\
 pK_s & \mbox{ otherwise }  \end{array} \right.,\; B,B'\in {\cal B},\end{equation}
Here $p = -1$ for the type 1 setting.

%\begin{equation} \label{Cmatrices} C_{BB';G}  =  \left\{ \begin{array}{ll}
%(s-1)K_s & \mbox{ if }  B' = B, \\
 %-K_s & \mbox{ otherwise }  \end{array} \right.,\; B,B' \in {\cal B}. \end{equation}

  For a fixed $B$, we eliminate all $\widehat{\alpha^{B'}}, B' \neq B$ from (\ref {ith}) by using (\ref {Kmatrix}) and (\ref
  {CmatricesIn(b)}). Then we get an equation involving only $\widehat{\alpha^{B}}$ and $\widehat{\alpha^V}$. We use this equation to eliminate all
  $\widehat{\alpha^{B}}$'s from (4.12). Then we get the  reduced normal equation for $\widehat{\alpha^V}$, which is $C_d \widehat{\alpha^V} =Q$, where
  $C_d$ is as in the statement.

  %{\bf Proof of (b):} In this case, the individual C-matrices are as given below.

%\begin{equation}\label{CmatricesIn(b)} C_{BB';G}  =  \left\{ \begin{array}{ll}
%(p+s)K_s & \mbox{ if }  B' = B, \\
 %pK_s & \mbox{ otherwise }  \end{array} \right.,\; B,B'\in {\cal B},\end{equation}

% Proceeding in the same way as in the proof of (a), we eliminate all $\widehat{\alpha^B}$ and
%get the reduced normal equation in $ \hat{\tau}$, the C-matrix of which is as in the statement.

{\bf Proof of (b):} Let $\tilde{{\cal B}} = {\cal B} \setminus \{E\}$.

 By the hypothesis,  the individual C-matrices are as follows.
\begin{eqnarray}\label{CmatricesIn(c)} C_{BB';G}  &=&  \left\{ \begin{array}{ll}
(s+1)K_s & \mbox{ if }  B' = B, \\
 K_s & \mbox{ otherwise }  \end{array} \right.,\; B,B' \in \tilde{{\cal B}},\\
 C_{BE;G} & = & \left \{ \begin{array}{ll}
sK_{s+1} & \mbox{ if }   B = E, \\
 0 & \mbox{ otherwise }  \end{array} \right. ,\; B \in {\cal B}.
 \end{eqnarray}

Following the same procedure as in Case (a)  we get the C-matrix as in the statement.  $\Box$

\section{Optimality of three multi-way designs}
\setcounter{equation}{0}

In this section we shall study a few multi-way designs, which are duals of main effect plans. So, we define what we mean by the dual of a
main effect plan.

\begin{defi}\label{dual} Consider a main effect plan ${\cal P}$ for a $s_1 \times \cdots s_m$ experiment on $n$ runs laid out on $b$
blocks. By the dual  $d ({\cal P})$ of ${\cal P}$ we mean  the following design in a multi-way heterogeneity set up. $d ({\cal P})$ has $b$
treatments,
to be tested on $n$ experimental units. The units are subjected to $m$ heterogeneity directions (blocking factors), the ith one having
$s_i$ levels.  A run $x = (x_1, \cdots x_n)'$ of  ${\cal P}$ corresponds to the  unit, say $x$, having $x_i$ as the level of the $i$th
blocking factor, $ 1 \leq i \leq m$. If $x$ is in the $j$th block of ${\cal P}$, then $d ({\cal P})$  allocates treatment $j$ is to $x$,
$1 \leq j \leq b$.
\end{defi}

%We assume that $s$ is a prime power.
We proceed to the description of the multi-way designs of interest.

\begin{nota} \label{setupDetails}  $s = ht +1$ is a prime power. $F$ will denote the field of order $s$. $F^*$ will denote the set of
non-zero elements of $F$.

(a)   For an integer $h$, $I_h =\{0,1,\cdots h-1\}$.

Throughout this section addition in the suffix will be  modulo $h$, whenever the suffix is in $I_h$.

(b) $C_i, i \in I_h$ will denote the cosets  of the subgroup $C_0$ of order $t$ of $F^*$.
  The cosets are so ordered  that $C_iC_j = C_{i+j}, \; i,j \in I_h$.
  %(addition in the suffix being modulo $h$).
  $\bar{C_i} = C_i \cup \{0\}$.

(c) $S_V = \{(x,i), x \in F, i \in I_h\}$ will be the set of treatments for all the designs we discuss below. Thus, $v = hs$.

%(c) Consider an $m$-way design $\triangle$ such that the set of levels of every blocking factor is $F$ or $F^+$ and the treatment set
%is $F \times I_h$. For $p \in F$, let $\tilde{p}$ denote the $m+1 \times 1$ vector $(p1'_m, (p,0))'$. Then the $m$-way design consisting of
%$s|\triangle|$ allotted units $\{\delta + \tilde{p} : \delta \in \triangle, \; p \in F \}$ will be denoted by $\triangle \oplus F$. WE may
%refer to this new design as the design ``generated from $\triangle$ by adding $F$".

%set up with $U$ as the set of all units. Then, $U \oplus F$ will denote the set of $s|U|$ units $\{u + a1_m, \;u \in U,
%a \in F\}$. Here addition in each co-ordinate is the field addition. This may be referred to as the set up `generated from $U$ by adding $F$'.

(d)  All factors except $E$ in Setting of type 3 have $F$ as the set of levels, while for $E$ it is $F^+ =
F \cup \{\infty\}$.

\end{nota}

\subsection{Optimality of $d^*_1$ :}  Corollary 2.3 of Morgan and Uddin (1996) has constructed   an universally optimal main effect plan,
say  $\rho_1$, for a $s^h$ experiment on $hs$ blocks of size $t$ each (here $s,h$ and $t$ are as in Notation \ref {setupDetails}). The dual
of $\rho_1$ is an $h$-way design. In this section we shall study this multi-way design, named $d^*_1$, regarding its performance in terms
of optimality. We present a description of $d^*_1$ in terms of our notation. Before that we describe a class of matrices.

 \begin{nota}\label{impMatrix}
  For $i \in I_h, L_i$ will denote the following $s \times s$ matrix with rows and columns indexed by $F$.
   $$L_i(x,y) = \left\{ \begin{array}{ll}
1 & \mbox{ if }  y - x \in  C_i, \\
 0 & \mbox{ otherwise } \end{array} \right.$$

 $L$ will denote the following partitioned matrix of order $v \times v$.
 $$L = ((L_{i-j}))_{i,j \in I_h}. $$
 %\mbox{(The addition in the suffix is modulo $h$)}.$$
 \end{nota}

{\bf The set up of  $d^*_1$:} The set of units is $U = F^* \times F$. The set ${\cal B}$ of block factors  is of size $h$; the members of
${\cal B}$ are as follows. Fix $P = \{p_i : i \in I_h\} \subset F^*$, such that $p_i \in C_i, i \in I_h$. ${\cal B} = P$.
The levels of $p_i$'s in different units are as follows.
\begin{equation} \label{levelBi*} \mbox{ the level of $p_i$ on the unit }(a,b) \in  F^* \times F \mbox{ is given by } \eta_2((a,b),p_i)
= ap_i + b.\end{equation}
Let us compute the $p_i \times p_j$ incidence matrix $M_{ij},\; i\neq j$ in the above set up. Fix $x,y \in F$. From (\ref {levelBi*}) we
see that the only possible unit $u = (a,b)$ such that $\eta_2(u,p_i) = x$ and $\eta_2(u,p_j) = y$ is given by
$$  a = \frac{x-y}{p_i - p_j}, b = \frac{p_i y - p_j x}{p_i - p_j}.$$
This is a unit only if $x \neq y$. Thus, $M_{ij} = J_s - I_s$, for  $i \neq j \in I_h$ and therefore  this is a setting of type 1 with $m =
h$.

 {\bf The design $d^*_1$ :}  We describe the allocation of treatments  to
 the units of the above set up [recall Notation \ref {setupDetails} (c)].
 \begin{equation} \label{treatmentd*1}\mbox{For } u = (a,b) \in U = F^* \times F, \; \; \eta_1(u) = (b,i), \mbox{ if } a \in C_i.
 \end{equation}
    %The unit $(a,b) \in  F^* \times F$ receives treatment $\eta_1 (a,b) = (b,j)$ if $a \in C_j$.
    [Recall that any $a \in F^*$ is in a unique $C_j, \; j \in I_h$].
    \vskip5pt
{\bf Computation of the treatment-versus-$p_i$ incidence matrix $N_i$ of $d^*_1$ :} Fix $i \in h, (x,j) \in F \times I_h, \; y \in F$.
$N_i ((x,j),y)$ is the number of $u \in  F^* \times F$ satisfying
\begin{equation} \label{ni(x,y)}\eta_1(u) = (x,j), \;\;\eta_2 (u,p_i) = y.\end{equation}
By (\ref {levelBi*}) and (\ref {treatmentd*1}), the only possible unit $u = (a,b)$ satisfying (\ref {ni(x,y)}) is given by
$$ a = (y-x)p^{-1}_i,\; b = x, \mbox{ provided }(y-x)p^{-1}_i \in C_j.$$ So,
$N_i ((x,j),y) = 1 $ if $ y-x \in C_{i+j}$ and  $0$ otherwise.
%    [The addition in the suffix is modulo $h$.]
 \begin{equation} \label{incMatd*1}  \mbox{Thus, } N_i =  \left [ \begin{array}{c} L_i \\ L_{i+1} \\ \vdots \\
  L_{i+h-1} \end{array} \right ],  \mbox{where $L_i$'s are as in Notation \ref {impMatrix}}. \end{equation}

%[Recall Notation \ref l{allotment}].

We see that $d^*_1$ is equireplicate with replication number $r = t$.

\begin{lem}\label{prop*d1}  $d^*_1$ satisfies  the following properties.

 (a) $C_{d^*_1} = rK_v - (1/s) LL' -  1/(s(s-h)) J_h \otimes K_s  + (hr^2/s^2) J_v. $

(b) The spectrum of $ C_{d^*_1}$ is $r^{h-1} (r - \frac{1}{s-h})^{s-1} (r-1)^{(h-1)(s-1)}0^1$.
\end{lem}

{\bf Proof :}  The proof of (a) follows from Lemma \ref {C-matTreatment} (a) in view of (\ref {incMatd*1}).  Proof of (b) is in the
Appendix.

Let $d$ be an equireplicate competing design in the setting of  $d^*_1$. We define $H = H_d$ by
\begin{equation} \label{H} H = \sum \limits_{B \in {\cal B}} N_B.\end{equation}

\begin{defi} An $m$-way design is said to be {\bf totally binary} if the  entries of $H$  are $0$ or $1$.
The class of all equireplicate and totally binary $m$-way designs in the setting of $d^*_1$ will be denoted by ${\cal D}^B_r$.
\end{defi}
We  note that any design in ${\cal D}^B_r$ shares the following property with $d^*_1$.

\begin{lem}\label{propH} For a design $d \in {\cal D}^B_r$, the following hold.

(a) $tr(C_d) = tr(C_{d^*_1})$.

(b)  By permuting the rows and columns if necessary, $H_d$ can be reduced to
$ H_d = 1_h \otimes ( J_s - I_s ).$

(c) The spectrum of $ S_{\cal B} S'_{\cal B}$ is $0^{v-s+1} \;h^{s-1}$. [Here $ S_{\cal B}$ is as in Lemma \ref {C-matTreatment} (a)]. \end{lem}

{\bf Proof :} (a) follows from the definition of ${\cal D}^B_r$. To prove (b) we see that
 for  a totally binary design, the entries of each $N_B$ are
$0$ or $1, \;\; B \in {\cal B}$, as the entries of each $N_B$ are non-negative.
Let $\bar{H} = J_{hs \times s} - H_d$. Since $d \in {\cal D}^B_r$, $\bar{H}$ is a $0,1$ matrix with exactly one entry $1$ in every row and
 exactly $h$ entries $1$ in every column. So, $\bar{H} = 1_h \otimes  I_s $ upto permutation of the rows and columns. Hence (b) follows.

 By the definition of $ S_{\cal B}$, the above relation imply that
 \begin{equation} \label{SB}
   S_{\cal B} S'_{\cal B} = HH' - ((s-1)^2/s) J_v  = J_h \otimes K_s. \end{equation}
  Now,  (c) is immediate. $\Box$

\vskip5pt

{\bf Remark 5.1:} Now onwards we shall assume  that for every deign $d$ in ${\cal D}^B_r$, the rows and columns of each $N_B$ is permuted
(if necessary) in accordance with the permutation (if any) used in Lemma \ref {propH}.

We proceed towards studying the performance of $d^*_1$.

\begin{theo}\label{opt d1*} $d_1^*$ is M-optimal in  ${\cal D}^B_r$, the class of all equireplicate and totally binary designs
 in the Setting of $d_1^*$.
\end{theo}

Towards the proof of this theorem, we fix a design  $d \in {\cal D}^B_r$  and study certain properties of the matrices in the expression
for $C_d$.
%(which is in statement (a) of Lemma \ref {C-matTreatment}). To begin with, we note that a setting containing $d$ satisfies
%$m = h$ and the replication number is $r = t$. We also note that

 We shall use the following compact notations. Let $N$ and $C_{\cal B}$ be the following $v \times sh$ matrices
$$N = [ N_B : B \in {\cal B} ],   \mbox{ and } C_{\cal B}  = [ C_B : B \in {\cal B}].$$
 The following relation hold.
\begin{equation}  \label{GiPi} C_{\cal B} C'_{\cal B} =   NN' - (hr^2/s)  J_v. \end{equation}
%($b) $1_v$ is an eigenvector of  $NN'$ and $HH'$ with eigenvalue $(s-1)^2$ and  $h(s-1)^2$ respectively.\end{lem}
\vskip5pt

%\begin{lem}\label{exprCd} For $d \in {\cal D}^B_r$, the C-matrix of $d$ is as follows.
%$$ C_d = r K_v - (1/s)NN' - 1/(s(s-h)) J_h \otimes K_s + (r(s-1)/s^2)J_v.$$ \end{lem}

%{\bf Proof :} For an equireplicate design, $C_0 = r K_v$. Again, in this set up,  $S_{\cal B} S'_{\cal B} = J_h \otimes K_s$
%by Lemma \ref {propH}. Moreover, $\sum \limits_{B \in {\cal B}} C_B C'_B =  C_{\cal B} C'_{\cal B} = NN' - (r(s-1)/s) J_v$. Hence the
%result follows from Lemma \ref {C-matTreatment} (a). $\Box$
%Let ${\cal N} (A)$ denote the null space of the matrix $A$ and $\nu (A)$ the nullity of $A$.

\begin{lem}\label{sum eValueP}  The sum of $s-1$ smallest positive eigen values of
 $C_{\cal B} C'_{\cal B}$ is  $ \leq s-1$. \end{lem}

 {\bf Proof : }% We recall that  $C_{\cal B} C'_{\cal B} =   NN' - (hr^2/s)  J_v$.
  Let $X = \{x_1, \cdots x_{s-1}\}$ be a set of orthonormal vectors in $ \langle 1_s  \rangle^{\perp}$. Let $Z =  \{ z_i = 1_h \otimes x_i, \; x_i \in X\}$.
 %Let $Z = 1_h \otimes x_{s \times 1} : x' 1_s = 0$.
   For a $z \in Z, \; Nz = Hx_i$ for some $i$, so that $Nz = -z$ and so  $z'N'Nz = z'z$. Let $\tilde{\mu}_1 \leq \cdots \leq
   \tilde{\mu}_{s-1}$ be the  smallest $s-1$ positive eigen values of $NN'$. Then, $\tilde{\mu}_1, \cdots , \tilde{\mu}_{s-1}$ are also the
   smallest $s-1$ positive eigen values of $N'N$. So,
 $$ \sum_{i=1}^{s-1} \tilde{\mu}_i \leq \sum_{i=1}^{s-1} (z'_i N'N z_i)/(z'_i z_i) = (s-1).$$

 Since  $N$ is an incidence matrix, the largest eigenvalue of $NN'$ corresponds to the eigenvector $1_v$. Thus, the eigenvector $e_i$
 corresponding to $\tilde{\mu}_i$ can not be $1_v$ and therefore $e'_i 1_v = 0$, for each $i$. Hence the result follows from (\ref {GiPi}).
  $\Box$

Let ${\cal N} (A)$ denote the null space of the matrix $A$ and $\nu (A)$ the nullity of $A$.

\begin{lem}\label{bigEvaluesN} $\mu(C_d)$ satisfies the following.

(0) $\mu_0(C_d) = 0$.

(i)  $\mu_{v-i}(C_d) = r, \; 1 \leq i \leq h-1$.

(ii)  $\sum_{i=0}^{s-2} \mu_{v-h-i} (C) \geq (s-1)(r - \frac{1}{s-h})$.
\end{lem}

{\bf Proof :}
%$$ C_d = r K_v - (1/s)NN' - 1/(s(s-h)) J_h \otimes K_s + (r(s-1)/s^2)J_v.$$
By Lemma \ref {C-matTreatment} (a) $$C_d = C_0 - \frac{1}{s} C_{\cal B} C'_{\cal B} - \frac{1}{s(s-h)} S_{{\cal B}} S'_{{\cal B}}. $$
Using  (\ref {CBcolsum0}) we see  that   $1'v C_{{\cal B}} = 0 = 1'_v S_{\cal B}$.
Thus, $C_d 1_v = 0$ and (0) is proved.

We now prove (i). Substituting for $C_{\cal B} C'_{\cal B}$ from (\ref {GiPi}) and $S_{\cal B} S'_{\cal B}$ from (\ref {SB})
  we get the following expression for $C_d$.
 %For that we use the fact that $C_{\cal B} C'_{\cal B} = NN' - (r(s-1)/s)J_v$.
$$ C_d = r K_v - (NN' - \frac{r(s-1)}{s}J_v) - \frac{1}{s(s-h)} (HH' - \frac{(s-1)^2}{s}J_v) .$$
 By definition of $H = H_d$, ${\cal N} (NN')  \subset {\cal N} (HH')$. Now,  let  $W = \langle 1_h \rangle^{\perp} \otimes 1_s$. Since $d$
 is equireplicate, $Nw = 0, \forall w \in W$. Therefore, $ \nu (NN') = \nu (N'N)  \geq  |W| = h-1$, implying $ \nu (HH') \geq h-1$.
 Let $ x \in {\cal N} (NN')$. Since  $1_v$ is the eigenvector of $NN'$ corresponding to the largest eigenvalue, $x \neq 1_v$,. Therefore, $x' 1_v = 0$. So, using the expression of $C_d$ above, (\ref {SB})and  (\ref {GiPi})  we see that $C_d x = r$. Hence (i) follows.

  Next we proceed to prove (ii), which is about the next largest eigen values of $C_d$,
 Let $P = a C_{\cal B} C'_{\cal B}   + b S_{\cal B} S'_{\cal B} , \; a,b >0$. While proving (i), we have also proved that $\mu_i (P) = 0,\;
 0 \leq i \leq h-1$. Let $\tilde{\mu}_1 (T) \leq \cdots \leq \tilde{\mu}_{s-1} (T)$ be the  smallest $s-1$ positive eigen values of $T,\; T
 = P$ or $C_{\cal B} C'_{\cal B}$. Fix $ i : 1 \leq i \leq s-1$.
 By an well-known result [see Corollary III.2.2 of Bhatia (2013), for instance] we get
 $$\tilde{\mu}_i (P) \leq a \tilde{\mu}_i (C_{\cal B} C'_{\cal B}) + b \mu_{v-1} (S_{{\cal B}} S'_{{\cal B}}),$$
 which is $a \tilde{\mu}_i (C_{\cal B} C'_{\cal B}) + bh$, by Lemma 5.2 (c). So, by Lemma 5.3,
  \begin{equation}  \label{sumSmallevalue} \sum_{i= 1}^{s-1} \tilde{\mu}_i (P) \leq (s-1) (a + bh).\end{equation}
  But $P$ becomes $r K_v - C_d$, if we put $a = \frac{1}{s}, b =  \frac{1}{s(s-h)}$. So, the result follows from (\ref {sumSmallevalue}).$\Box$

  \vskip5pt

{\bf Proof of Theorem \ref {opt d1*} :} From Lemma \ref {prop*d1} (b) the spectrum of $C_{d^*_1}$ is
$$r^{h-1} (r - \frac{1}{s-h})^{s-1} (r-1)^{(h-1)(s-1)}0^1.$$
Now let us  take a pair of $v - h \times 1$ vectors $x$ and $y$, where $y_i = \mu_i(C_d)$ and $x_i = \mu_i(C_{d^*_1}), \; 1 \leq i \leq
v-h$. Lemma 5.4  says that the other eigenvalues of $C_d$ are equal to the corresponding ones of  $C_{d^*_1}$. This,
in view of Lemma \ref {propH} (a)
, shows that  Lemma \ref {suffMopt} is applicable here. Applying that lemma and the results obtained in Lemma
5.4 to $x$ and  $y$ we find that $y$ is M-worse than $x$.  Hence the proof is complete. $\Box$

\subsection{Optimality of $d^*_2$} Among the  universally optimal main effect plans constructed in Lemma 2.5  of Morgan and Uddin (1996) one is a plan, say  $\rho_2$, for a $s^t$ experiment on $hs$ blocks of size $t+1$ each (here $s,h$ and $t$ are as in Notation \ref
{setupDetails}). We shall study  the mutiway design, named $d^*_2$, which is the dual of $\rho_2$
We describe it using our notation.

{\bf The setting :} The set of units $U = \{(a,b,j) : j \in I_h, b \in F, a \in \bar{C_j}\}$. The set of block factors is ${\cal B} =
C_0$. The set of level of each  block factor is $F$. Finally,
\begin{equation} \label{levelB2} \mbox{ the level of $\alpha \in C_0$ \; on the unit } u = (a,b,j) \in  U \mbox{ is given by } \eta_2(u,
\alpha) = a\alpha + b.\end{equation}
{\bf Computation of the $\alpha$ -versus- $\beta$ incidence matrix $M_{\alpha\beta}$:} Fix $\alpha \neq \beta, \; \alpha, \beta \in C_0$.
For $x,y \in F$, the only possible unit  $u = (a,b,i)$ satisfying
$\eta_2(u,\alpha) = x, \;\; \eta_2(u,\beta) = y$ is given by
 $$a = \frac{x-y}{\alpha- \beta}, \; b =  \frac{ay-bx}{\alpha- \beta} \mbox{ and }i \in I_h \mbox{ is such that }a \in \bar{C}_i. $$
 Therefore,
$ m_{\alpha, \beta} (x,y) = 1$  if   $x \neq y$ and $h$  if   $x = y$.  Thus ,
$$ M_{\alpha, \beta} = (h-1) I_s + J_s, \; \forall \alpha \neq \beta.$$
 Hence, this setting is of type 1 with $m = t$.

 {\bf The design $d^*_2$ :} $S_V = F \times I_h$. For $u = (a,b,j) \in U, \; \eta_1(u) = (b,j)$. Thus, for $\alpha \in C_0, \; (x,i) \in
 S_V$,
 the unit $u = (a,b,j)$ satisfies $\eta_1(u) = (x,i),\;  \eta_2(u,\alpha) = y$ only if $b = x, j = i, a = (y-x) \alpha^{-1}$ and $a$ is in
 $\bar{C}_i$.
 It follows that for every $\alpha \in C_0$, the $((x,i), y)$th entry of the treatment-versus-$\alpha$ incidence matrix $N_\alpha$ is
 $ N_{\alpha} ((x,i), y) = 1$  if  $ y - x \in  \bar{C}_i$ and $0$  otherwise. Thus,
 \begin{equation} \label{Nalpha} \mbox{ For each } \alpha \in C_o, \; N_\alpha =  \left [
 \begin{array}{cc} L_0 +I   \\ L_{1} + I \\ \vdots \\ L_{h-1} + I \end{array} \right ], \mbox{ where $L_i$ is as in Notation \ref
 {impMatrix}}.\end{equation}

We present the following well-known result for ready reference.

  \begin{lem}\label{spectrBIBD} If $N$ is the incidence matrix of a BIBD with parameters $(v,b,r,k,\lambda)$, then the spectrum of $N'N$ is
  is  $0^{b-v}(r -\lambda)^{v-1} (rk)^1$. \end{lem}

\begin{lem}\label{d*2prop}
(a) The C-matrix of $d^*_2$ is $C_{d^*_2} = rK_v - (t/u) (N'_2 N_2 - ((t+1)^2/s) J_{v \times s})$, where $N_2$ is the incidence matrix of a
BIBD with parameters $(v = s, b = hs, r = h(t+1), k = t+1, \lambda = t+1)$.

(b)  The spectrum of $ C_{d*_2}$ is $r^{(h-1)(s-1)} (r - (t/u)(t+1))^{s-1} 0^1$.
\end{lem}

 {\bf Proof  :} From the description of $d^*_2$, one can see that $d^*_2$ is a design in Setting 2  with $m = t, p = h-1$. So, (a)
  follows from Lemma \ref {C-matTreatment} (a) in view of (\ref {Nalpha})

  (b) follows from (a) in view of Lemma \ref {spectrBIBD}.  $\Box$

\vskip5pt

We now proceed towards studying optimality aspects of $d^*_2$. We shall use a few well-known results, presented below.

\begin{lem}\label{minSQ} Suppose $x_1, \cdots x_n$ are real numbers with $\sum_{i=1}^{n} x_{i} = T$. Then, the following hold.
(a)  $\sum_{i=1}^{n} x^2_{i} \geq T^2/n, ``=" \mbox{ when } x_i = T/n, \; \forall i$.

(b) In particular, if  $x_i$'s are integers, then  $\sum_{i=1}^{n} x^2_{i}$ is minimum, when $x_i = [T/n]$ or $[T/n] + 1$. \end{lem}

\begin{lem}\label{AAmin} Consider an $m \times n$ matrix $A$ having the $m \times 1$  vector  $b$ as the vector of row sums. Then,
$  AA' \geq bb'/n.$
\end{lem}

{\bf Proof :} Fix an arbitrary $x \in R^m$, let $y = A'x$. Now, using Lemma \ref {minSQ} on the components of $y$, we get the result.
$\Box$

\begin{lem}\label{sumDiffMat} Consider  $m \times n$ matrices $A_1, \cdots A_q$. Let  $T = \sum_{i=1}^q A_i$. Then,
the following hold.

(a) $ \sum_{i=1}^q A_i A'_i \geq (1/q) TT'.$

(b) $q \sum_{i=1}^{q} A_i A'_i - T T' = \sum_{i=1}^{q} \sum \limits_{ j < i} (A_i - A_j)  (A_i - A_j)'.$
\end{lem}

{\bf Proof  of (a) :} Fix an arbitrary $x \in R^m$. For $ 1 \leq i \leq q$, let us write $A'_ix$ as $y_i = \left [ \begin{array}{cccc}
y_{i1} & \cdots y_{in} \end{array}\right]^\prime$ . Then, $x' \sum_{i=1}^q A_i A'_i x =   \sum_{i=1}^q y'_i y_i = \sum_{j=1}^n \sum_{i=1}^q
y^2_{ij}$.
Now, $\sum_{i=1}^q y_{ij} = \sum_{l=1}^m T(l,j) x_l = z_j$, say. Therefore, by Lemma \ref {minSQ},  $\sum_{i=1}^q y^2_{ij} \geq z^2_j /q$.
Since $z_j$ is the $j$th entry of $T'x$, the result follows.

(b) Follows by straightforward computation. $\Box$

The next two results are not so well-known.
\begin{lem}\label{MinTrHH'} Suppose $ H = \sum_{l=1}^{m} N_l$, where each $N_l$ is an integer matrix of order $v \times b$. Suppose each
$N_l$ satisfies  $\sum_{j=1}^{b} N_l (i,j) = r,\; i = 1, \cdots v$, where $r < b$. Then,  $Tr(HH') \geq m^2 vr.$
\end{lem}

{\bf Proof :} $  Tr(HH') =  \sum_{i=1}^{v} \sum_{j=1}^{b} (H(i,j))^2.$ Again, $\sum_{i=1}^{v} \sum_{j=1}^{b} H(i,j) = mvr$.
 Therefore, by Lemma \ref {minSQ}, $\sum_{i=1}^{v} \sum_{j=1}^{b} (H(i,j))^2$ is minimum, when the following hold.
 $H (i,j) =  [rm/b] \mbox{ or } [rm/b] + 1,  1 \leq i \leq v, 1 \leq j \leq b.$
Since a sufficient condition for the above is
$ N_1 = N_2 = \cdots N_m, \mbox{ and } N_1 (i,j) = 0 \mbox{ or } 1,  \;\; 1 \leq i \leq v, 1 \leq j \leq b$,  the result follows. $\Box$

\vskip5pt

% The next result is an important step  towards proving the optimality properties of $d^*_2$ and $d^*_4$.

\begin{lem}\label{gammaMbetter2}  Let $A$ be an  $n \times n$ n.n.d matrices each with row sum $0$, rank $\leq \rho$ and trace $\geq T$.
Let $C = dK_n - A$, where $d \geq T/\rho = a$ (say) and $K_n$ is as in (\ref {Kmatrix}). Let $\gamma \in (I\!\! R^+)^n$ be the
vector given by $\gamma_0 = 0, \gamma_i = d - a$ for $1 \leq i \leq \rho, \; \; \gamma_i = d $ for $\rho + 1 \leq i \leq n-1$.
 Then, $\gamma$ is M-better than $\mu(C).$ \end{lem}

{\bf Proof :} By definition of $C, \; \mu_0(C)  = 0 = \gamma_0, \; \;  \mu_i (C) = d -  \mu_{n+1-i} (A)$ for $i > 0$.
Therefore, $\sum_{i=1}^{\rho} \mu_{i} (C) = d \rho - tr(A) \leq \rho (d - a)$.
Since $(1/l)\sum_{i=1}^l \mu_{i} (C)$ is increasing in $l$, it follows that
$\sum_{i=1}^l \mu_{i} (C) \leq l(d - a), \; 1 \leq l \leq \rho$.
Since $\mu_i(C) \leq d \; \forall i$, the result follows. $\Box$

\begin{theo} \label{d2*} $d_2^*$ is  M-optimal in  ${\cal D}_2$, the class of all equireplicate designs in the setting containing $d_2^*$.
\end{theo}

{\bf Proof :} Let $d \in {\cal D}_2$. Then, $C_d$ is  as given in (a) of  Lemma \ref {C-matTreatment} with $m = t, p = h-1$.
 So, using Lemma \ref {sumDiffMat} (b) on the set of matrices $\{C_B, B \in B \in {\cal B} \}$,  we get
 $$ C_d = C_0 - (1/u) \sum \limits_{B \in {\cal B}} C_BC'_B - (h-1)/(su) \sum \limits_{B' \neq B} D_{BB'}, \; u = s + (h-1)t, $$
where $D_{BB'} = (C_B-C_{B'}) (C_B-C_{B'})', \; B' \neq B \in {\cal B}$. Since $d$ is equireplicate, $C_0 = rK_v$. Now,
applying   Lemma \ref {sumDiffMat} (a) on the set of matrices $ \{ D_{BB'}, B,B' \in {\cal B}\}$  we get $C_d \leq C_1 = rK_v - (1/u)
\sum \limits_{B \in {\cal B}} C_B$.  Again, by  Lemma \ref {sumDiffMat} (a) on the set of matrices $ \{ C_{B}, B \in {\cal B}\}$  we get
$C_1 \leq C_2$,   where  $C_2  = rK_v - (1/(tu) S_{{\cal B}}S'_{{\cal B}}$.
Here $S_{{\cal B}}$ is as in Lemma \ref {C-matTreatment} (a). Now, $C_2$ is of the form of $C$ in Lemma \ref {gammaMbetter2} with $n = hs,
d = r, A = (1/tu)S_{{\cal B}}S'_{{\cal B}}, \rho = s-1$. Also,  $T = (h-1)(s-1)(t+1)t/u$ by Lemma \ref {MinTrHH'}. In view of (\ref
{CBcolsum0}) the row sum of $A$ is $0$. So, $a = (h-1)(t+1)t/u$. By  Lemma \ref {gammaMbetter2}, $\mu(C_2)$ is M-worse than $\gamma$. But
from Lemma \ref {d*2prop} we see that $\gamma = \mu(C_{d^*_2})$. Since $C_d \leq C_2$, the result follows from Lemma \ref {suffMopt0}.
$\Box$

\subsection{Construction and Optimality of $d^*_3$} In this section we construct a multi-way design $d^*_3$ and prove its optimality
property.
In this section we assume that $s \equiv 3 \pmod 4$. We follow Notation \ref {setupDetails} with the extra assumption that $h = 2$. So,
here $C_0$ (respectively $C_1$) is the set of non-zero squares  (respectively  non-zero non-squares). Moreover, $t = (s-1)/2$.

The following result lies at the foundation  of the construction of $d^*_3$.
\begin{lem}\label{factorsInW} Let $s \equiv 3 \pmod 4$ be a prime power. Then there is a subset $W$ of $C_0$ and a function $f : W
\rightarrow C_1$ satisfying the following.

(a) $|W|  = (s-3)/4$.

(b) For every $\xi \in W$, $ (\xi - 1) (f(\xi) - 1) \in C_0$.

(c) For $\xi \neq \xi' \in W$,  $(\xi - \xi') ( f(\xi)- f(\xi')) \in C_0.$
 \end{lem}

{\bf Proof :} Let $W = \{x \in C_0 : 1-x^2 \in C_0\}, \;  \tilde{W} = \{x \in C_1 : 1-x^2 \in C_1\}$. Note that $x \rightarrow -x$ is a
bijection from $W$ onto $(C_1 \setminus \tilde{W})\setminus \{-1\}$. Therefore, $|W| = |C_1 \setminus \tilde{W}| -1 = (s-3)/2 -
|\tilde{W}|$.
Thus, $|W| + |\tilde{W}|  = (s-3)/2$. Also, $x \rightarrow -1/x$ is a bijection from $W$ onto $\tilde{W}$. Thus, $|W| = |\tilde{W}|$, which
 proves (a).

Let   $f : W \rightarrow  \tilde{W}$ be defined by $f(x) = -1/x,\; x \in W$. Then, for $\xi \in W$, $1-\xi^2 \in C_0$, so that
$(1 -\xi) (1 - f(\xi)) = \xi^{-1}  (1 - \xi^2)$ which is in $C_0$. This proves (b).

Again, for $\xi \neq \xi' \in W$, $(f(\xi) - f(\xi'))(\xi - \xi') = (\xi \xi')^{-1} (\xi - \xi')^2 \in C_0$, which implies (c). $\Box$

\vspace{.5em}

We shall now construct the design  $d^*_3$. Let $w = |W|$. So, $w = (t-1)/2$.

{\bf The set up :} $U = \{(a,b,i) : a \in \bar{C_0}, b \in F, i = 0,1\}$.  ${\cal B} = W \cup \{\infty\}$, where $W$ is as in Lemma \ref
{factorsInW}. The set of levels of factor $\infty$ is $F^+= F \cup \{\infty\}$, while the  set of levels of the factor $\xi $ is $F$, for
each $\xi \in W$. The  levels of different factors in the units are as follows. [Here $f$ is as in the proof of Lemma \ref {factorsInW}].

  For $u = (a,b,i) \in U,\; \xi \in W$,
\begin{eqnarray} \label{allotTr xi}
 \eta_2 (u,\xi) = &\left\{ \begin{array}{ll}
 b + a\xi   & \mbox{ if }  i = 0,  \\
b - a f(\xi)  & \mbox{ if }  i = 1 \end{array} \right. \mbox{ and }\\
 \eta_2 (u,\infty) = &\left\{ \begin{array}{ll}
a + b  & \mbox{ if }  a \in \bar{C_0}, \; i = 0,  \\
b -  a  & \mbox{ if } a \in C_0,\; i = 1,   \\
 \infty & \mbox{ if } a = 0, \; i = 1  \end{array} \right. \end{eqnarray}
%Let $w = |W$, where $W$ is an in Lemma \ref {factorsInW}.

%\begin{lem}\label{type3} The set up described above is   a setting of type 3 with $m = w+1$.
%\end{lem} {\bf Proof :}
Fix $\xi \in W$.  For $x \in F, \; y \in F^+$. We count the number of $u = (a,b,i) \in U$ such that
\begin{equation}  \label{M xiinfty} \eta_2(u,\xi) = x \mbox{ and } \eta_2 (u,\infty) = y. \end{equation}
If $y = \infty$, then by the equation next to (\ref {allotTr xi}), $u = (0,x,1)$ is the only unit satisfying  (\ref  {M xiinfty}).
%Thus, $m_{\xi\infty} (x,\infty) = 1, x \in F$.

 Now, let $ y \in F$. If $y = x$, then  $u = (0,x,0)$ is the only unit satisfying (\ref  {M xiinfty}).

So, let $ y \neq x, x,y \in F$.  Using (\ref {allotTr xi})
and the equation next to it we find  that the only unit $u$ satisfying (\ref  {M xiinfty}) is
$$ u = \left\{ \begin{array}{ll} (\frac{y-x} {1 - \xi}, \frac{x - y \xi}{1 - \xi}, 0) & \mbox{ if } (1 - \xi)(y-x) \in C_0 \\
(\frac{y-x} {f(\xi) -1},\frac{ y f(\xi) -x}{f(\xi) -1}, 1) &  \mbox{if } (f(\xi) -1)(y-x) \in C_0. \end{array} \right. $$
 Moreover, by Lemma \ref {factorsInW} (b), exactly one of $1 - \xi$ and $f(\xi) - 1$ is in $C_0$. Thus for a given pair $x \neq y$,
  exactly one of the two cases occur.
   %Thus,  the first entry of $u$ is in $C_0$, as is required.
   Hence $M_{\xi \infty} = J_{s \times s+1}$.

Now, fix $\xi' \neq \xi \in W$.  We fix a pair $x,y \in F$ and count the number of $u = (a,b,i) \in U$ such that
\begin{equation}  \label{M xixi'} \eta_2(u,\xi) = x \mbox{ and } \eta_2 (u,\xi') = y. \end{equation}

If $y = x$, then $u_1 = (0,x,0)$ and $u_2 = (0,x,1)$ are the only units satisfying (\ref {M xixi'}). Let
 $y \neq x$.  From (\ref {allotTr xi}) we find that the only unit $u$ satisfying (\ref  {M xixi'}) is
$$  u = \left\{ \begin{array}{ll} (\frac{x - y} {\xi - \xi'}, \frac{y \xi - x \xi'}{\xi - \xi'}, 0)  & \mbox{ if } (\xi - \xi')(x-y)\in C_0
\\
 (\frac{x - y} {f(\xi') - f(\xi)}, \frac{y f(\xi) - x f(\xi')} {f(\xi') - f(\xi)}, 1)  & \mbox{ if}  (f(\xi') - f(\xi))(x-y)\in C_0.
 \end{array} \right. $$
 By by Lemma \ref {factorsInW} (c) exactly one of $\xi - \xi'$ and $f(\xi') - f(\xi)$ is in $C_0$, so that for any given pair $x \neq y$,
 exactly one of these two cases occur.

   Hence  $M_{\xi\xi'} = I_s + J_s$ and
   the setting is  a setting of type 3 with $m = (s+1)/4 = (t+1)/2$.

{\bf The design  $d^*_3$ :} $d^*_3$ has $v = 2s$ treatments ; $S_V = \{(x,i) : x \in F,\; i = 0,1\}$. The treatments are assigned to the
units by the rule that for $u = (a,b,i)$
\begin{equation}  \label{trAllot}\eta_1(u) = (b,i).\end{equation}

\begin{theo}\label{propOfd^*3}
%Suppose $s \equiv 3 \pmod 4$. Then there exists a design
  $d^*_3$ is  a design in a setting of Type 3  with $(s+1)/4$ blocking factors, satisfying the following properties.

(a) $d^*_3$ is  equireplicate with replication number $r = t+1 = (s+1)/2$.

(b)  The treatment-versus-blocking factor  incidence matrices are as follows. [The set of blocking factors is indexed by $W \cup
\{\infty\}$ , where $W$ is as in  Lemma \ref {factorsInW}]

(i) For all $\xi \in W, N'_{\xi}$ is the incidence matrix of a BIBD, (independent of $\xi$) with parameters $(v = s, b = 2s, r = s+1, k =
t+1, \lambda = t+1)$.

(ii) $N'_{\infty}$ is the incidence matrix of a BIBD with parameters  $(v = s+1, b = 2s, r = s, k = t+1, \lambda = t)$.

(c)  For every $\xi \in W, \; \xi$ and $\infty$ are adjusted orthogonal with respect to the treatment factor.
%[Recall Definition \ref {adjustedOrth}]

 (d) The spectrum of $ C_{d^*_3}$ is $(\frac{rs}{s+w})^{s-1} (\frac{r(s-1)}{s})^s 0^1$.
% $(r - wr/u)^{s-1} (r-1)^s 0^1$.
\end{theo}

{\bf Proof :} By the description of $d^*_3$, we see that there are $w = (s-3)/4$ $s$-level and one $s+1$-level blocking factors.
(a) follows from the allocation of treatments [see (\ref {trAllot})]. We proceed to prove (b). Fix $\xi \in W$.  For a fixed $(x,i) \in S_V$ and $y \in F$, $N_{\xi} ((x,i),y))$ is the number of units $u$ such that
\begin{equation}  \label{N xicomp}  \eta_1(u) = (x,i) \mbox{ and }  \eta_2(u,\xi) = y. \end{equation}
Using (\ref {allotTr xi}) and (\ref {trAllot}) we see that (\ref {N xicomp} ) has a solution only when $y - x \in \bar{C_0}$. In that case
the only unit $u$ satisfying (\ref {N xicomp} ) is
$$ u = \left\{ \begin{array}{ll} ((y-x)/\xi,x,0)  \mbox{ if } i = 0 \\
((x-y)/f(\xi),x,1)  \mbox{ if } i = 1.  \end{array} \right.$$
Since $\xi$ and $-f(\xi)$ are in $C_0$, the first entry of $u$ is in $\bar{C_0}$, as required.
Therefore,
\begin{equation} \label{N xi}
N_\xi =  \left [ \begin{array}{l} L_0 + I  \\
 L_0 + I    \end{array} \right ],\; \xi \in W. \end{equation}
 where $L_0$ is as in Notation \ref {impMatrix} with $h = 2$. Thus, (b) (i) is proved.

 Towards proving (b) (ii), we fix $(x,i) \in S_V$ and $y \in F^+$.
$N_{\infty} ((x,i),y)$ is the number of units $u$ such that
\begin{equation}  \label{Ninftycomp} \eta_1(u) = (x,i) \mbox{ and }  \eta_2(u,\infty) = y. \end{equation}
Clearly if $y = \infty$, then $u = (0,x,1)$ is the only unit satisfying  (\ref {Ninftycomp}) when $i = 1$ and there is no such unit when $i
= 0$.

Next, let $y \in F$. $y = x \implies$ that the only unit    satisfying (\ref {Ninftycomp}) is $u = (0,x,0)$. Let $ y \neq x$.
From the equation next to (\ref {allotTr xi}) and (\ref {trAllot}) we see that  if $y = x$, then  the only unit satisfying (\ref
{Ninftycomp}) is
$$ u = \left\{ \begin{array}{ll} (y-x,x,0)  \mbox{ if }  y - x \in C_0 \\
  (x-y, x, 1) \mbox{ if }  y - x \in C_1. \end{array} \right.$$
  It follows that
 \begin{equation} \label{Ninfty} N_{\infty} =  \left [ \begin{array}{ll} L_0 + I & 0_{s \times 1} \\
 L_1     & 1_s  \end{array} \right ], \end{equation}
 where $L_0$ and $L_1$ are as in Notation \ref {impMatrix} with $h = 2$. Hence (b) (ii) is proved.

 Now we prove (c). By definition of $L_0$ and $L_1$ given in  Notation \ref {impMatrix},  $(L_0 + I)1_s = (t + 1) 1_s$ and $L_0 + I + L_1 = J_s$. So,  from (\ref  {N xi}) and (\ref {Ninfty}) we see that$N'_{\xi} N_{\infty} = (L_0 + I)' J_{s \times s+1} = (t + 1) J_{s \times s+1} = r M_{\xi \infty}$ for every
 $\xi \in W$. Now (c) follows from Definition 2.3.

Finally, we prove (d). By Lemma \ref  {C-matTreatment} (b), Lemma \ref {sumDiffMat} (b)
and the  information on $d^*_3$ we have got so far, we see that
   \begin{equation}  \label{Cd*3} C_{d*_3} = r K_v - \frac{w}{s+w} C_{\xi}  C'_{\xi} -  \frac{1}{s} C_{\infty} C'_{\infty},\end{equation}
    where $C_B$ is as in (\ref {CBcolsum0}), $B \in  W \cup \{\infty\}$. So,  $ C_{\xi}  C'_{\xi} = N_{\xi}N_{\xi}' - \frac{r^2}{s}J_v $ and $C_{\infty} C'_{\infty} = N_{\infty}(N_{\infty})' - \frac{r^2}{s+1}J_v$. Here  $N_{\xi}$ is as in (\ref {N xi}) and $N_{\infty}$ is as in (\ref {Ninfty}). By  part (b) and Lemma \ref {spectrBIBD} we find that the spectrum of $C_{\xi}  C'_{\xi}$ is $r^{s-1} 0^{v-s+1}$ and  the spectrum of $C_{\infty}  C'_{\infty}$ is $r^s 0^{v-s}$.  Now (d) follows from (\ref {Cd*3}), in view of Lemma \ref {adjOrth}. $\Box$

 \begin{theo} \label{optd3*}  $d_3^*$ is  M-optimal in  ${\cal D}_3$, the class of all equireplicate designs in the setting of $d_3^*$.
\end{theo}

{\bf Proof :} Let $d \in {\cal D}_3$. From Lemma \ref {C-matTreatment} (b) we get $C_d$ is  as given below.
 $$C_d = C_0  -  \frac{1}{s}\sum \limits_{\xi \in  W^*} C_{\xi} C'_{\xi}
+ \frac{1}{su} \tilde{S}_W \tilde{S}^{\prime}_W,  $$

where $u = s + w, \; W^* =  W \cup \{\infty\}$ and $ \tilde{S}_W$ is defined in the same way as in Lemma \ref {C-matTreatment} (b).
 %\end{document}

% $$ C_d = C_0 - (1/u) \sum \limits_{B \in {\cal B}} C_BC'_B - (h-1)/(su) \sum \limits_{B' \neq B} D_{BB'}, \; u = s + (h-1)t, $$
%where $D_{BB'} = (C_B-C_{B'}) (C_B-C_{B'})', \; B' \neq B, B,B'\in {\cal B}$.

By arguments similar to those used in Theorem \ref {d2*} we see that $C_d \leq C_1$, where
$$C_1 = rK_v  - (1/u) \sum \limits_{\xi \in  W} C_{\xi} C'_{\xi}  -  (1/s) C_\infty C'_{\infty}.$$
Again by similar arguments we find that $C_1 \leq C_2$, where
$$C_2= rK_v -  (1/wu)\tilde{S}_W \tilde{S}'_W   .$$

  We see that $C_2$ is of the form of $C$ in Lemma \ref {gammaMbetter2} with $n = 2s, A = (wu)^{-1} \tilde{S}_W\tilde{S}'_W, d = r, \rho =
  s-1, T = wr(s-1)/u.$ So, $a = wr/u$. Therefore, by the same lemma, $\mu (C_2)$ is M-worse than $\gamma$. As $C_1 \leq C_2, \: \mu (C_1)$ is
    M-worse than $\gamma$. In particular,
 $$\sum_{j=1}^{l} \mu_j (C_1) \leq l (r  - a), 1 \leq l \leq s-1. $$
Again, one can check that $Tr(C_1) \leq (s-1) (r - a)  + r(s-1)$.
%where $\delta_1 = r - \alpha $ and $\delta_2 = (r-1)/s$.
 Let $\delta$ be a $v-1 \times 1$ vector such that
$\delta_i =  r - a$ if $ 1 \leq i \leq s-1$  and  $r(s-1)/s$  if  $s \leq i \leq v-1$.
It follows that $\mu (C_1)$ is M-worse than $\delta$.
But, by Theorem  \ref {propOfd^*3} (d) $\delta = \mu (C_{d^*_3} )$.
 Since $C_d \leq C_1$  the result follows. $\Box$

\vspace{.5em}

{\bf Remark 5.2:} Theorem \ref {optd3*} is an extension of the result of Bagchi and Shah (1989). The design analogous to $d^*_3$ in a
row-column set up exists for every prime power $s$. [See  Preece, Wallis and  Yucas (2005) and Nilson and Cameron (2017) for more details].

\vspace{.5em}
{\bf Remark 5.3:} In Bagchi and Bagchi (2020) an MEP  ${\cal P}^*$ is constructed for an $s^t(s+1)$ experiment on $2s$ blocks of size $t+1$
each. Let $d_4$ denote the dual of $\rho^*$, which is a multi-way design with $t+1$ blocking factors, one having $s+1$ levels and others
having $s$ levels. Although $d_4$  satisfies  properties similar to those of $d^*_3$ given in Theorem \ref {propOfd^*3}, it does not satisfy condition (c). As a result, we are unable to prove M-optimality of $d_4$. Whether $d_4$ satisfies any specific optimality property remains to be seen.
%\end{document}

\section{Appendix}

\setcounter{equation}{0}
Consider $C_{d^*_1}$ of Section 5.1. We prove the following result.
\begin{theo}\label{spectrN(1)} The spectrum of  $C_{d^*_1}$ is $r^{h-1} (r - (1/(s-h))^{s-1} (r-1)^{(h-1)(s-1)}.$ \end{theo}

 To prove this,  we need a number of tools.

\begin{nota} (0) $s = p^m = ht+1$, where $p$ is a prime, $m,t,h$ are integers, $m \geq 1, h,t \geq 2$. $F_p$ and $F_s$ are finite fields of
orders $p$ and $s$ respectively. [In this section we use the notation $F_s$ (rather than $F$ like in the other sections) so as to
distinguish it from the field of order $p$].

(i) Addition and subtraction in $I_h = \{0.1, \cdots h-1\}$ will always be modulo $h$.

The rows and columns of every $s \times s$  (respectively $h \times h$) matrix will be indexed by $F_s$ (respectively $I_h$).
Moreover, the rows and columns of every $hs \times hs$   matrix will be indexed by $I_h \otimes X $.

(ii) $\eta$ and $\omega$ are primitive $h$th and $p$th roots of unity.

(iii) Consider the function trace  $ : F_s \rightarrow F_p$  defined as follows. $trace(x) = \sum_{i=1}^{m} x^{p^i},\; x \in F_s$.
[This is $F_p$-linear and into $F_p$ since $x \rightarrow x^p$ is an automorphism of $F_s$ and its fixed field is $F_p$].

(iv) $U$ and $V$ are unitary matrices of orders $h$ and $s$ respectively, given as follows.
$$U(i,j) = (1/\sqrt{h}) \eta^{ij}, i,j \in I_h \mbox{ and } V(x,y) = (1/\sqrt{s}) \omega^{trace(xy)}, x,y \in F_s. $$

(iv) Consider the  sums $g_i$ given by $g_i = \sum \limits_{ x \in C_i} \omega^{trace(x)}, i \in I_h$.

(v) As in Notation \ref {impMatrix}, $L$ will denote the $hs \times hs$ matrix $(( L_{i-j}))_{i,j \in I_h}$, where for
 $i \in I_h, L_i$ will denote the $s \times s$ matrix
   $L_i(x,y) = \left\{ \begin{array}{ll}
1 & \mbox{ if }  y - x \in  C_i, \\0 & \mbox{ otherwise } \end{array} \right.$

(vi) For $k \in I_h$, $G_k$ is the $h \times h$ matrix given by $G_k(i,j) = g_{i-j+k}, \; i,j \in I_h$.

(vii)  For $k \in I_h$, $E_k$ is the $s \times s$ diagonal matrix given by
$$E_k(x,x) = \left \{ \begin{array}{ll} - t   & \mbox{ if } x =0 \\
1 & \mbox{ if } x \in C_k \\
0 & \mbox{ otherwise } \end{array} \right.$$
 Here $C_k$ is as in Notation \ref {setupDetails}.
 \end{nota}

We study the behaviour of $L$ under the actions of $U$ and $V$.

\begin{lem}\label{VNV} $(I_h \otimes V)^* L (I_h \otimes V) = \sum \limits_{ k \in I_h} G_k \otimes E_k$.
\end{lem}
{\bf Proof :} By definition of $L$, the left hand side of the statement is
$((V^* M_{i-j} V))_{i,j \in I_h}$.
%(( C_{ij;U}))_{i,j \in T}

{\bf Claim !} $V^* M_i V = Diag(\lambda_i(x),\; x \in F_s)  = \sum \limits_{ k \in I_h} g_{i+k} E_k.$
Here $ \lambda_i (x) = \left \{ \begin{array}{ll}  & \mbox{ if } x =0 \\
 g_{i_k} & \mbox{ if } x \in C_k.\end{array} \right.$

The equality can be verified by computation.  Now,
$$\sum \limits_{ k \in I_h} g_{i+k} = \sum \limits_{ k \in I_h} g_k = \sum \limits_{ x \in F^*_s} \omega^{trace(x)} =-1
\mbox{ as } \sum \limits_{ x \in F_s} \omega^{trace(x)} = 0. $$

So, by definition of $E_k$, ${\cal V}_2 (1,1) = t$ and for $x \in C_j, {\cal V}_2 (x,x) = \sum \limits_{ k \in I_h} g_{i+k} \delta_{jk} =
g_{i+j}$. Thus,  the claim is proved and the second equality follows.

Now, the result follows from the definition of $G_k$.

\begin{nota} (a) $W$ will denote the $hs \times hs$ unitary matrix $U \otimes V$.

(b) $T$ is the $h \times h$ diagonal matrix with the  entries :
$T(l,l) = \eta^l, \; l \in I_h$.
\end{nota}

 \begin{lem}\label{spectrN} $W^* L W$ is a diagonal matrix with the following entries. For $i \in I_h, x \in F_s$, the $(i,x)$th diagonal
 entry of  $W^* L W$
 is $$\delta (i,x) =  \left \{ \begin{array}{ll} s-1   & \mbox{ if } i=0, x =0 \\
0 & \mbox{ if } i \neq 0, x =0 \\
\eta^{ik}\sum \limits_{ j \in I_h} g_j \eta^{-ij} & \mbox{ if } x \in C_k. \end{array} \right.$$
\end{lem}

{\bf Proof :} Since $W = (I \otimes v) (U \otimes I) $, by Lemma \ref {VNV},  $W^* L W = (U \otimes I)^* (\sum \limits_{ k \in I_h}
G_k \otimes E_k )(U \otimes I) = \sum \limits_{ k \in I_h} (U^* G_k U) \otimes E_k$.  One can verify that
$$U^* G_k U =  \sum \limits_{ j \in I_h} g_{k-j} T^j, \; k  \in I_h.$$
So, $W^* L W = \sum \limits_{ k \in I_h} \sum \limits_{ j \in I_h}  g_{k-j} T^j \otimes E_k$. Now, the formulae for $T$ and $E_k$ imply the
result. $\Box$

We are to obtain $W^* L (L)' W = (W^* L W)^2$. So,  we need to find a simpler expression for $|\sum \limits_{ j \in I_h} \eta^{-ij}|^2$.
For that we need the following.
\begin{nota} (a) $\Omega_h$ is the multiplicative group of all $h$th roots of unity.

(b) For $i \in I_h$, $\chi_i :  F^*_s \rightarrow \Omega_h$ is defined by $\chi_i (x) = \eta^{-ij}$, if $x \in C_j, j \in I_h$.
\end{nota}

\begin{lem}\label{gjEtaij2} $|\sum \limits_{ j \in I_h} g_i \eta^{-ij}|^2 =  \left \{ \begin{array}{ll} 1   & \mbox{ if } i=0,\\
  s  & \mbox{ if } 0 < i < h \end{array} \right.$ \end{lem}

 {\bf Proof :} Since $C_jC_k = C_{j+k}$ (where the addition in the suffix is modulo $h$), $\chi_i$ is a group homomorphism (character) on
 $F^*_s$.  From the definition of $g_j$'s we have
 $$ \sum \limits_{ j \in I_h} g_j \eta^{-ij} = \sum \limits_{j \in I_h} \sum \limits_{x \in C_j} \omega^{trace(x)} \chi_i(x) =  \sum
 \limits_{x \in F_s} \omega^{trace(x)}  = g(\chi_i),$$
 which is the Gauss sum attached to the character $\chi_i$. But $|g(\chi_i)|^2 = \left \{ \begin{array}{ll} 1   & \mbox{ if } i=0,\\
  s  & \mbox{ if } i \neq 0 \end{array} \right.$, by a classical result on such Gauss sums [see, for instance Chapter 10 of Ireland and
  Rosen (1982)]. Hence the result. $\Box$

  Putting the information from Lemmas \ref {spectrN} and \ref {gjEtaij2} together, we get the spectrum of $L (L)'$.

  \begin{lem}\label{spectrumNN'}
   $W^* L (L)' W= D$, where the diagonal entries of the diagonal matrix $D$ are as follows.
  $$ |\delta (i,x)|^2 =  \left \{ \begin{array}{ll} (s-1)^2   & \mbox{ if } i=0, x =0 \\
0 & \mbox{ if } i \neq 0, x =0 \\
 1   & \mbox{ if } i=0, x \neq 0, \\
  s  & \mbox{ if } i \neq 0, x \neq 0. \end{array} \right.$$
    \end{lem}

 In order to get the  spectrum of $C_{d_1^*}$ we need the spectrum of $HH'$.

 \begin{lem}\label{spectrumHH'}  $W^* HH' W$ is a diagonal matrix with the $(i, x)$th entry $= \left \{ \begin{array}{ll} h(s-1)^2   &
 \mbox{ if } i=0, x =0 \\
  h & \mbox{ if } i \neq 0, x \neq 0 \\
 0   & \mbox{ if } i \neq 0. \end{array} \right.$
  \end{lem}

  {\bf Proof :} Let $\triangle_h$ denote the $h \times h$ matrix having the $(0,0)$th entry $1$ and all other entries $0$. $\triangle_s$
  is defined in a similar manner. It is easy to verify that
  $$ U^* J_h U = h \triangle_h \mbox{ and } V^* J_s V = s \triangle_s.$$

 Since $ H = 1_h \otimes ( J_s - I_s )$, [recall  Lemma \ref {propH}], we see that
 $$W^* HH' W =   h \triangle_h \otimes (s \triangle_s - I_s)^2,$$
 which is a diagonal matrix with the entries as in the statement. $\Box$

 {\bf Proof of Theorem \ref {spectrN(1)} :}
 Lemmas \ref {spectrumNN'} and \ref {spectrumHH'} imply the result in view of the expression for $C_{d^*_1}$ in Lemma \ref {prop*d1}.
 $\Box$

 \section{Reference} \begin{enumerate}

\item Bagchi, B. and Bagchi, S. (2001). Optimality of partial geometric designs. Ann. Stat., vol. 29, p : 577 - 594.

\item Bagchi, S. and Bagchi, B.  (2020). Aspects of optimality of plans orthogonal through other factors. Submitted.

\item  Bagchi, S. and Mukhopadhyay, A.C. (1989) Optimality in the presence of two-factor interactions among the nuisance factors. Comm.
    Stat. Theory. Method. vol. 18, p : 1139-1152.

\item  Bagchi, S. and Shah, K.R. (1989). On the optimality of a class of row-column designs. Jour. Stat. Plan. Inf., vol. 23, p : 397-402.

\item Bhatia. R. (2013) Matrix Analysis. Graduate texts in Mathematics, Springer.

%\item Bose, M. and  Bagchi, S. (2007). Optimal main effect plans in blocks of small size, Jour. Stat. Prob. Let., vol.  77, p : 142-147.

%\item Das, A. and Dey, A. (2004). Optimal main effect plans with nonorthogonal blocks. Sankhya, vol. 66, p : 378-384.

%\item Dey, A. and Mukherjee, R. (1999). Fractional factorial plans. Wiley Series in probability and Statistics.

%\item Eccleston, J.A. and Russell, K.G. (1975). Connectedness and orthogonality in multi-factor designs. Biometrika,  vol. 62, p :341-345.

  \item Eccleston, J.A. and Russell, K.G. (1977). Adjusted orthogonality in nonorthogonal designs. Biometrika,  vol. 64, p : 339-345.

\item Ireland, K. and Rosen, M (1982). A classical introduction to modern number theory, Springer Verlag.

\item Kiefer, J. (1975). Construction and optimality of generalized Youden designs. In: Srivastava JN (ed)
A survey of statistical design and linear models. North-Holland, Amsterdam, p : 333 - 353

\item Marshall, A.W, Olkin, I, Arnold, B.C. (2011) Inequalities: theory of majorization and its applications.
Springer series in statistics, 2nd edn. Springer, New York.

\item Morgan, J.P. (1997).  Optimal design for interacting blocks with OAVS incidence. Metrika, 45, p : 67-83.

\item Morgan, J.P. and Uddin, N. (1996). Optimal blocked main effect
plans with nested rows and columns and related designs. Ann. Stat.
vol. 24,   p:  1185-1208.

\item Mukerjee, R., Dey, A. and Chatterjee, K. (2001). Optimal main effect plans with non-orthogonal blocking. Biometrika, 89,  p:
    225-229.

\item Mukhopadhyay, A.C. and  Mukhopadhyay, S. (1984). Optimality in a balanced multi-way heterogeneity set up. In Proceedings of the
    Indian Statistical Institute Golden Jubilee International conference on Statistics : Applications and new directions., p : 466-477.

\item Nilson, T. and   Cameron, P.J. (2017) Triple arrays from difference sets. Jour. Combin Designs. 2017, 25, p : 494 - 506.

\item Preece,  D. A., Wallis,  W. D. and  Yucas, J. L. (2005). Paley triple arrays.   Australas J. Combin., 33, p: 237 - 246.

\item Shah, K.R. and Eccleston, J.A. (1986). On some aspects of row-column designs.  Jour. Stat. Plan. Inf. 15, p : 87-95.

\item Shah KR, Sinha BK (1989) Theory of optimal designs.Lecture notes in statistics, vol 54. Springer. Berlin

\end{enumerate}

* Foot note : Both the authors are retired from Indian Statistical Institute, Bangalore Center

\vskip5pt

On behalf of all authors, the corresponding author states that there is no conflict of interest.
\end{document}